\input amstex \documentstyle{amsppt}\magnification=1200\nologo
\NoBlackBoxes  \hsize=16 truecm \vsize = 9 truein \voffset = 7.5
truemm \hoffset = 0.1 truemm  \topmatter

\title
A speciality theorem for curves in $\bold P^5$
\endtitle

\author Vincenzo Di Gennaro* and Davide Franco** \endauthor

\leftheadtext{Vincenzo Di Gennaro and Davide Franco}
\rightheadtext{A speciality theorem for curves in $\bold P^5$}

\medskip
\abstract \nofrills

ABSTRACT. Let $C\subset \bold P^r$ be an integral projective
curve. One defines the speciality index $e(C)$ of $C$ as the
maximal integer $t$ such that $h^0(C,\omega_C(-t))>0$, where
$\omega_C$ denotes the dualizing sheaf of $C$. Extending a
classical result of Halphen concerning the speciality of a space
curve, in the present paper we prove that if $C\subset \bold P^5$
is an integral degree $d$ curve not contained in any surface of
degree $< s$, in any threefold of degree $<t$, and in any fourfold
of degree $<u$, and if $d>>s>>t>>u\geq 1$, then $ e(C)\leq
{\frac{d}{s}}+{\frac{s}{t}}+{\frac{t}{u}}+u-6. $ Moreover equality
holds if and only if $C$ is a complete intersection of
hypersurfaces of degrees $u$, ${\frac{t}{u}}$, ${\frac{s}{t}}$ and
${\frac{d}{s}}$. We give also some partial results in the general
case $C\subset \bold P^r$, $r\geq 3$.



\endabstract
\endtopmatter

{\footnote""{*Dipartimento di Matematica, Universit\`a\ di Roma
\lq \lq Tor Vergata\rq\rq, 00133 Roma, Italy. E-mail:
digennar{\@}axp.mat.uniroma2.it}

{\footnote""{ **Dipartimento di Matematica e Applicazioni \lq\lq
R. Caccioppoli\rq\rq, Universit\`a\ di Napoli \lq \lq Federico
II\rq\rq, Ple Tecchio 80, 80125 Napoli, Italy. E-mail:
dfranco{\@}unina.it}}

\bigskip

Let $C\subset \bold P^r$ be an integral projective curve. One
defines the speciality index $e(C)$ of $C$ as the maximal integer
$t$ such that $h^0(C,\omega_C(-t))>0$, where $\omega_C$ denotes
the dualizing sheaf of $C$. In [GP] Gruson and Peskine prove the
following theorem concerning the speciality index of an integral
space curve:
\bigskip
\proclaim{Theorem A} Let $C\subset \bold P^3$ be an integral
degree $d$ curve not contained in any surface of degree $< s$.
Then one has:
$$
e(C)\leq {\frac{d}{s}}+s-4.
$$
Moreover equality holds if and only if $C$ is a complete
intersection of surfaces of degrees $s$ and ${\frac{d}{s}}$.
\endproclaim \bigskip
In the present paper we extend this theorem to curves in $\bold
P^5$, in the following sense:
\bigskip
\proclaim{Theorem B} Let $C\subset \bold P^5$ be an integral
degree $d$ curve not contained in any surface of degree $< s$, in
any threefold of degree $<t$, and in any fourfold of degree $<u$.
Assume $d>>s>>t>>u\geq 1$. Then one has:
$$
e(C)\leq {\frac{d}{s}}+{\frac{s}{t}}+{\frac{t}{u}}+u-6.
$$
Moreover equality holds if and only if $C$ is a complete
intersection of hypersurfaces of degrees $u$, ${\frac{t}{u}}$,
${\frac{s}{t}}$ and ${\frac{d}{s}}$.
\endproclaim \bigskip

In Remark $(i)$ below,  we make explicit  the condition
$d>>s>>t>>u$.

Taking into account the general inequality:
$$
deg(C)\cdot e(C)\leq 2p_a(C)-2 \tag 1
$$
($p_a(C)=$ arithmetic genus of $C$), in [GP], Remarque $3.6$, the
authors notice that Theorem A is an immediate consequence of the
Halphen bound for the genus of space curves and the consequent
classification of the curves with maximal genus.  Similarly, as we
will see in the proof, in the case $C\subset \bold P^4$, i.e. when
$u=1$, our Theorem B  is an easy consequence of the analogue of
the Halphen bound for the genus of curves in $\bold P^4$ verifying
{\it{flag conditions}}, shown in [CCD].

Given integers $s_1,s_2,\dots,s_{l}$, with $1\leq l\leq r-1$ and
$s_l\geq r-l+1$, we say that a nondegenerate integral projective
curve $C\subset \bold P^r$ verifies the flag condition
$(r;s_1,s_2,\dots,s_{l})$ if $C$ has degree $s_1$, and, for any
$i=2,\dots,l$, it is not contained in any integral, projective
subvariety of $\bold P^r$ of dimension $i$ and degree $<s_i$.
Except for the case $l\leq 2$ (see [CCD3]), for the case $l=3$ and
$r=4$ (see [CCD]), and for some partial results in the case $l\geq
3$ (see [CCD2]), it is not known an analogue of the Halphen bound
for curves verifying flag conditions in $\bold P^r$, $r\geq 5$.
Therefore one cannot obtain a speciality theorem in $\bold P^r$,
$r\geq 5$, as a direct consequence of  Halphen bound as indicated
before.

In the case $r=5$ we are able to overcome this difficulty in the
following way. First we notice that, under the hypotheses of
Theorem B, from $(1)$ it follows that if $e(C)\geq
{\frac{d}{s}}+{\frac{s}{t}}+{\frac{t}{u}}+u-6$, then the
arithmetic genus of $C$ is \lq\lq sufficiently big\rq \rq to apply
[CCD2], Theorem $2.2$, and so to deduce that $C$ belongs to a flag
such as
$$
C\subset S\subset T\subset U\subset {\bold P^5}, \tag 2
$$
where $S$, $T$ and $U$ are a surface, a threefold and a fourfold
of degrees $s$, $t$ and $u$. Next, using a sort of a {\it{
numerical coarse adjunction formula}} on $S$ (see Lemma $1$
below), we are able to compare the sectional genus of $S$ with the
Halphen bound for curves in $\bold P^4$ quoted before, and to
deduce that $S$ is in fact a complete intersection (for the proof
of Lemma $1$, which relies on Castelnuovo Theory, we refer to [D];
see also Remark $(ii)$ below for further comments). At this point
Theorem B follows from Lemma $2$ below, which is a generalized
\lq\lq weak\rq\rq version of Theorem A, for curves contained in a
complete intersection surface. The proof of Lemma $2$ relies on a
liaison argument, and does not need the numerical assumption on
the degrees.

Since the proof of Theorem B depends on  results obtained in
[CCD], [CCD2] and [D], we are forced to assume that $r=5$, and
that $d>>s>>t>>u$. We use this last condition also in some
numerical estimates: see $(8)$, $(9)$ and $(13)$ below.

In the general case concerning curves in ${\bold P^r}$, $r\geq 3$
and $1\leq l\leq r-1$,  we notice the following  partial result.

\proclaim{Proposition 1} Let $C\subset \bold P^r$ be a reduced and
irreducible curve verifying the flag condition
$(r;s_1,\dots,s_{l})$. Assume that $s_1>>s_2>>\dots >>s_{l}$. Then
one has:
\smallskip
$(i)$
$$
e(C)\leq {\frac{s_1}{s_2}}+{\frac{s_2}{s_3}}+\dots +
{\frac{s_{l}}{s_{l+1}}}-l-{\frac{2}{r-l}}+{\frac{3}{4}};
$$
in particular when $l=r-1$ then
$$
e(C)\leq {\frac{s_1}{s_2}}+{\frac{s_2}{s_3}}+\dots +
{\frac{s_{r-2}}{s_{r-1}}}+s_{r-1}-(r+1)+{\frac{3}{4}};
$$
\smallskip
$(ii)$ if $l=r-1$ and $s_i$ divides $s_{i-1}$ for any
$i=2,\dots,r-1$, then
$$
e(C)\leq {\frac{s_1}{s_2}}+{\frac{s_2}{s_3}}+\dots +
{\frac{s_{r-2}}{s_{r-1}}}+s_{r-1}-(r+1)
$$
and, in such case, the bound is sharp;
\smallskip
$(iii)$
$$
e(C)\leq
{\frac{s_1}{s_2}}+{\frac{2G(r-1;s_2,\dots,s_l)-2-s_2}{s_2}}, \tag
3
$$
where $G(r-1;s_2,\dots,s_l)$ denotes the maximal arithmetic genus
for curves verifying the flag condition $(r-1;s_2,\dots,s_{l})$.
Moreover, if there exists some  subcanonical curve verifying the
flag condition $(r-1;s_2,\dots,s_{l})$ and with maximal genus
$G(r-1;s_2,\dots,s_l)$, then the bound is sharp.
\endproclaim \bigskip

When $l=2$, $G(r-1;s_2)$ denotes the Castelnuovo bound for
nondegenerate curves in $\bold P^{r-1}$ of degree $s_2$ (see also
Proposition $2$ below). In Remark $(i)$ below, we make explicit
the condition $s_1>>s_2>>\dots >>s_{r-1}$. Recall that a
projective curve $C$ is {\it{subcanonical}} if $\omega_C=\Cal
O_{C}(e(C))$. For integral curves, this is equivalent to require
that the equality holds in $(1)$. For instance, a Castelnuovo
curve in $\bold P^{r-1}$  of degree $s_2\equiv 2$ $mod(r-2)$ is
subcanonical. This gives an example for $(iii)$. See also
Proposition $2$ below for further examples for flag conditions of
length $\leq 3$, and Remark $(iv)$ below for a reformulation of
property $(iii)$.

We will see that Proposition $1$  is an easy consequence of
results contained in [CCD2] and [D]. Notice that the bound $(ii)$
is sharp because it is attained by complete intersection curves of
type $(s_{r-1}, {\frac{s_{r-2}}{s_{r-1}}},\dots,
{\frac{s_1}{s_2}})$. One expects that the bound $(ii)$ holds true
in any case, i.e. without the assumption $s_1>>s_2>>\dots
>>s_{r-1}\geq 2$, and that $s_i$ divides $s_{i-1}$ for any
$i=2,\dots,r-1$. Also, one expects that the bound $(ii)$ it is
attained {\it{only}} by  complete intersection curves. We are not
able to prove a  speciality theorem in this generality, neither to
eliminate the positive term ${\frac{3}{4}}$ appearing in $(i)$ of
Proposition $1$ (however we note that, by a suitable choice of the
numerical assumptions on the $s_i$'s, one may substitute the term
${\frac{3}{4}}$ with any {\it{positive}} number (compare with
[CCD2], Remark $1.3$)). Notice also that property $(i)$ of
Proposition $1$ improves the bound in [CCD2], Proposition $2.4$,
$ii)$.

Finally we recall that, extending previous results of Eisenbud and
Harris [EH], in [CCD3] one generalizes Halphen bound also for the
genus of curves in $\bold P^r$ of given degree $d$ and  not
contained in surfaces of degree less than a fixed number $s$.
Consequently, as we said before for Theorem A, by $(1)$ one easily
obtains a corresponding bound for the speciality. Therefore, in
the case of flag conditions of length $2$ (i.e. $l=2$), one has  a
more precise result than previous Proposition $1$. In fact, fix
integers $r$, $d$ and $s$, with $s\geq r-1\geq 2$, and define $m$,
$\epsilon$, $w$ and $v$ by dividing $d-1=ms+\epsilon$, $0\leq
\epsilon \leq s-1$, and $s-1=w(r-2)+v$, $0\leq v \leq r-2$. Denote
by $G(r;d,s)$ the maximal arithmetic genus for an integral
nondegenerate projective curve $C\subset \bold P^r$ of degree $d$,
not contained in any surface of degree $<s$. One may prove that,
when $d>>s$, then
$$
G(r;d,s)={\frac{d^2}{2s}}+{\frac{d}{2s}}\left[2G(r-1;s)-2-s\right]+R,
\tag 4
$$
where $G(r-1;s)={{w}\choose{2}}(r-2)+wv$ is the Castelnuovo bound
for a nondegenerate curve of degree $s$ in $\bold P^{r-1}$, and
$$
R={\frac{1+\epsilon}{2s}}(s+1-\epsilon-2G(r-1;s))+w(\epsilon-\delta)
-k{{w+1}\choose{2}}+{{\delta}\choose{2}}.
$$
The numbers $k$ and $\delta$ are defined by dividing
$\epsilon=kw+\delta$, $0\leq \delta<w$, when $\epsilon <w(r-1-v)$.
Otherwise define $k$ and $\delta$ by dividing
$\epsilon+r-2-v=(w+1)k+\delta$, $0\leq \delta<w+1$ (compare with
[CCD3]). Taking into account $(1)$ and the previous expression for
$G(r;d,s)$, with an elementary numerical analysis based on [CCD3]
we are able to prove the following:
\bigskip
\proclaim{Proposition 2} Fix integers $r$, $d$ and $s$, with
$s\geq r-1\geq 2$, and $d>>s$. Let $C\subset \bold P^r$ be an
integral nondegenerate projective curve of degree $d$, not
contained  in any surface of degree $<s$. Define $m$, $\epsilon$,
$w$ and $v$ by dividing $d-1=ms+\epsilon$, $0\leq \epsilon \leq
s-1$, and $s-1=w(r-2)+v$, $0\leq v \leq r-2$. Put
$G(r-1;s)={{w}\choose{2}}(r-2)+wv$. Then one has

$$
e(C)\leq {\frac{d}{s}}+{\frac{2G(r-1;s)-2-s}{s}}. \tag 5
$$
\medskip
Moreover the following properties are equivalent:
\medskip
$(i)$ $e(C)= {\frac{d}{s}}+{\frac{2G(r-1;s)-2-s}{s}}$;
\medskip
$(ii)$ $p_a(C)= G(r;d,s)$, and $v=0$ and $\epsilon=w$, or $v\geq
1$ and $\epsilon=w(r-1-v)+1$;
\medskip
$(iii)$ $p_a(C)= G(r;d,s)$ and $C$ is subcanonical.
\endproclaim \bigskip

As before, in Remark $(i)$ below we make explicit the condition
$d>>s$.

Now we are going to prove the announced results. We begin with the
proof of Theorem B.

\bigskip

To this purpose, first we examine the case $u\geq 2$, i.e. we
assume that $C\subset \bold P^5$ is nondegenerate.

Assume also that
$$
e(C)\geq {\frac{d}{s}}+{\frac{s}{t}}+{\frac{t}{u}}+u-6. \tag 6
$$
Then from $(1)$ we deduce that
$$
p_a(C) \geq
{\frac{d^2}{2s}}+{\frac{d}{2}}\left[{\frac{s}{t}}+{\frac{t}{u}}+u-6\right]+1.
\tag 7
$$
On the other hand, since $C$ satisfies the flag condition
$(5;d,s,t,u)$ and $d>>s>>t>>u$, from [CCD2], Proposition $1.1$ and
Theorem $2.2$, we know that
$$
p_a(C)\leq G^h(5;d,s,t,u):=
{\frac{d^2}{2s}}+{\frac{d}{2}}\left[{\frac{s}{t}}+{\frac{t}{u}}+u-6+\eta\right]
+\rho+1,
$$
where $\eta$ and $\rho$ are suitable rational numbers  such that
$$
|\eta|\leq {\frac{3}{4}} \quad {\text{and}} \quad |\rho|\leq
33{\frac{s^3}{t^2}}.
$$
From $(7)$ and our hypotheses $d>>s>>t>>u$, we get
$$
G^h(5;d,s,t,u)-p_a(C) \leq {\frac{d}{2}}\eta+\rho <
{\frac{d}{2}}\left[{\frac{t}{u(u+1)}} -3u\right]. \tag 8
$$
Therefore we may apply  [CCD2], Proposition $1.4$ and Theorem
$2.2$, ii), and deduce the existence of a flag such as $(2)$.

Our next step is to prove that $S$ and $T$ are in fact complete
intersections. To this aim, we are going to analyze the sectional
arithmetic genus $\pi$ of $S$ and to compare it with the quoted
bound for the genus of curves verifying the flag condition
$(4;s,t,u)$ in $\bold P^4$, appearing in [CCD].

First  we need the following lemma. For its proof we refer to [D],
Lemma.
\bigskip
\proclaim{Lemma 1}  Let $S\subset \bold P^r$ be an integral,
nondegenerate projective surface, of degree $s\geq r-1\geq 2$.
Denote by $\pi$ the sectional arithmetic genus of $S$. Let
$C\subset S$ be an integral, nondegenerate projective curve, of
degree $d\geq s^2$. Denote by $p_a(C)$ the arithmetic genus of
$C$. Then one has:
$$
 p_a(C)\leq
{\frac{d^2}{2s}}+{\frac{d}{2s}}(2\pi-2-s)+s^3/(r-2).
$$
In particular, when $d>2s^4/r-2$, one has
$$
e(C)\leq {\frac{d}{s}}+{\frac{2\pi-2-s}{s}}. \tag 9
$$
\endproclaim \bigskip

\noindent Comparing $(9)$ with $(6)$ we deduce that
$$
\pi \geq
{\frac{s^2}{2t}}+{\frac{s}{2}}\left[{\frac{t}{u}}+u-5\right]+1.
\tag 10
$$
On the other hand, by [CCD] we know that
$$
\pi \leq G(4;s,t,u), \tag 11
$$
where
$$
G(4;s,t,u):=
{\frac{s^2}{2t}}+{\frac{s}{2}}\left[{\frac{t}{u}}+u-5-{\frac{(u-1-\beta)(1+\beta)(u-1)}{ut}}\right]
+\rho_1(s,t,u)+1.
$$
In the previous expression the number $\beta$ is defined by
dividing $t-1=\alpha u+\beta$, $0\leq \beta <u$. The number
$\rho_1(s,t,u)$ is a rational number which only depends on $s$,
$t$ and $u$. For its definition we refer to [CCD]. For the moment
we only need  that
$$
|\rho_1(s,t,u)+1|\leq t^3/3. \tag 12
$$
Previous estimate follows from Lemma $1$, taking into account that
one may interpret  $G(4;s,t,u)$ as the arithmetic genus attained
by curves in $\bold P^4$ lying in surfaces of degree $t$, whose
sectional genus is equal to
${\frac{t^2}{2u}}+{\frac{t}{2}}(u-4)-{\frac{(u-1-\beta)(1+\beta)(u-1)}{2u}}+1$,
i.e. is equal to the Halphen bound for curves of degree $t$, not
contained in surfaces of degree $<u$ in $\bold P^3$ (see [CCD]).
From $(10)$, $(11)$ and $(12)$, and our assumption $s>>t$, we
obtain that
$$
u-1-\beta=0 \tag 13
$$
i.e. $u$ divides $t$, and that $\rho_1(s,t,u)\geq 0$. On the other
hand, from [CD], pg.2708-2709, we know that, when $u$ divides $t$,
then $\rho_1(s,t,u)\leq 0$, and $\rho_1(s,t,u)=0$ if and only if
$t$ divides $s$. It follows that $u$ divides $t$, $t$ divides $s$
and, from $(10)$ and $(11)$, that
$$
\pi =
{\frac{s^2}{2t}}+{\frac{s}{2}}\left[{\frac{t}{u}}+u-5\right]+1=G(4;s,t,u),
$$
i.e. the general hyperplane section $S'$ of $S$ has maximal
arithmetic genus $G(4;s,t,u)$ with respect to the flag condition
$(4;s,t,u)$. This in turn implies, from [CCD] again, that $T$ is a
complete intersection on $U$. Moreover also $S$ is a  complete
intersection on $T$. In fact, from the Hilbert function of $S'$,
we know that $S'$ is contained in a space surface  of degree
${\frac{s}{t}}$ not containing the general hyperplane section of
$T$, and such a surface  lifts to $\bold P^4$ because $S'$ is
arithmetically Cohen-Macaulay (see [CCD], Theorem and pg. 132).

At this point, taking into account $(6)$, to conclude the proof of
Theorem B in the case $u\geq 2$, it suffices to prove the
following lemma. We will prove it below.

\bigskip
\proclaim{Lemma 2}  Let $T\subset \bold P^r$ be an irreducible,
reduced, projective threefold of degree $t$. Assume that $T$ is
arithmetically Cohen-Macaulay and subcanonical, i.e.
$\omega_T=\Cal O_T(k)$ for some $k\in\bold Z$. Let $S\subset T$ be
an integral surface of degree $s$, complete intersection on $T$.
Let $C\subset S$ be an integral curve of degree $d$. Then one has:
$$
e(C)\leq {\frac{d}{s}}+{\frac{s}{t}}+k.
$$
Moreover equality holds if and only if $C$ is a complete
intersection of $S$ with a hypersurface of degree ${\frac{d}{s}}$.
\endproclaim \bigskip

This concludes the proof of Theorem B when $u\geq 2$. When $u=1$,
taking into account [GP], we may assume that $C$ is nondegenerate
in $\bold P^4$. In this case the proof is easier because it does
not need either Lemma $1$ or Lemma $2$. In fact, from [CCD], we
know that $p_a(C)\leq G(4;d,s,t)$ (compare with $(11)$ and
$(12)$). Therefore, assuming $(6)$ and using $(1)$, the same
previous analysis of the number $\rho_1(d,s,t)$ shows that $t$
divides $s$, $s$ divides $d$ and that $C$ is a complete
intersection. This concludes the proof of Theorem B.
\bigskip
Now we are going to prove Lemma $2$.

To this purpose, consider a hypersurface $P\subset \bold P^r$ of
degree $p>>0$ containing $C$, denote by $D$ the complete
intersection $S\cap P$, and denote by $R$ the residual scheme of
$C$ with respect to $D$. Since $T$ is subcanonical, the local
Noether sequence of liaison corresponding to $C$ and $R$ gives:
$$
0\to \Cal I_D\to \Cal I_R\to \omega_C(-k-{\frac{s}{t}}-p)\to 0
$$
($\Cal I =$ ideal sheaf). Since $D$ is arithmetically
Cohen-Macaulay and $h^0(C,\omega_C(-e(C))>0$, from previous
sequence it follows that there exists a hypersurface $Q$ of degree
$q=k+{\frac{s}{t}}+p-e(C)$ containing $R$ and not containing $D$.
The  complete intersection curve $D'=S\cap Q$ contains $R$. Set
$D'=R\cup R'$ the corresponding, possibly algebraic, linkage.
Since $deg(R')=d-s(e(C)-k-{\frac{s}{t}})$ and  $deg(R')\geq 0$, we
obtain the claimed inequality. If the equality holds, then $R'$ is
empty and $R=D'$ is a complete intersection, with $\omega_R=\Cal O
_R(k+{\frac{s}{t}}+q)$. Coming back to the first linkage, we find
$$
0\to \Cal I_D\to \Cal I_C\to \Cal O_R(q-p)\to 0.
$$
Therefore there exists a hypersurface  of degree $p-q$ containing
$C$ and not containing $D$. Since $d=s(p-q)$, then $C$ is the
complete intersection of $S$ with such a hypersurface. This
concludes the proof of Lemma $2$.
\bigskip \bigskip

As for the proof of Proposition $1$, we first notice that
properties $(i)$ and $(ii)$ simply follow form $(1)$ and from
[CCD2], Proposition $1.1$, $(v)$, and Theorem $2.2$. It remains to
prove property $(iii)$. To this purpose, denote by
$G(r;s_1,\dots,s_l)$ the maximal arithmetic genus among the curves
verifying the flag condition $(r;s_1,\dots,s_l)$. Then, from
$(1)$, we have
$$
e(C)\leq {\frac{2G(r;s_1,\dots,s_l)-2}{s_1}}.
$$
From [D], Theorem, $(d)$, and the fact that $s_1>>s_2$, we deduce
the inequality $(3)$. In order to prove the sharpness, let
$E\subset \bold P^{r-1}$ be a subcanonical curve verifying the
flag condition $(r-1;s_2,\dots,s_l)$ and of maximal arithmetic
genus $p_a(E)=G(r-1;s_2,\dots,s_l)$. We have
$e(E)=(2G(r-1;s_2,\dots,s_l)-2)/s_2$. By [D], Theorem, $(a)$, we
know that $E$ is arithmetically Cohen-Macaulay. Let $C$ be a
complete intersection of  the cone over $E$ in $\bold P^{r}$ with
a general hypersurface of degree $a>>0$. The curve $C$ is
arithmetically Cohen-Macaulay and, by [D], Theorem, $(b)$, it
verifies the flag condition $(r;s_1,\dots,s_l)$, where $s_1=as_2$.
We have to prove that $e(C)=a+e(E)-1$. By $(3)$, we already know
that $e(C)\leq a+e(E)-1$. On the other hand, since $E$ is
arithmetically Cohen-Macaulay, then  from the exact sequence:
$$
0\to \Cal O_E(e(E))\to \Cal O_E(a+e(E)-1)\to \Cal O_{\Gamma}\to 0,
$$
one deduces that $h_{\Gamma}(a+e(E)+1)=s_1-1$ ($h_\Gamma$ =
Hilbert function of the general hyperplane section $\Gamma$ of
$C$) . It follows that $e(C)\geq a+e(E)-1$: in fact, since $C$ is
arithmetically Cohen-Macaulay, then its speciality index is equal
to $max\{i: h_\Gamma(i)<s_1\}-1$. This concludes the proof of
Proposition $1$.
\bigskip \bigskip

Now we are going to prove Proposition $2$.

With the same notation as in $(4)$, a direct computation proves
that $|R|\leq s^3/(r-2)$. This also follows from Lemma $1$,
because $G(r;d,s)$ is attained by curves lying in Castelnuovo
surfaces of degree $s$. Therefore, since $d>>s$ and $p_a(C)\leq
G(r;d,s)$,  from $(1)$ we deduce the bound $(5)$.

Another way to obtain $(5)$ is the following: denote by
$h_\Gamma(i)$ the Hilbert function of the general hyperplane
section $\Gamma$ of $C$, and by $h(i)$ the Hilbert function of the
general hyperplane section of any nondegenerate curve of degree
$d$, not contained in surfaces of degree $<s$, with maximal genus
$G(r;d,s)$. In [CCD3] one proves that this function $h(i)$ depends
only on $r$, $d$ and $s$, and it may be explicitly computed. Since
the speciality index of $C$ is $\leq max\{i: h_\Gamma(i)<d\}-1$,
and $h_\Gamma(i)\geq h(i)$ for any $i$, it follows that $e(C)\leq
max\{i: h(i)<d\}-1$. An elementary computation proves that this
last number is bounded as in $(5)$.

Now, if the equality holds in $(5)$, then the number
$e(C)={\frac{d}{s}}+{\frac{2G(r-1;s)-2-s}{s}}$ is an integer. A
computation proves that this is equivalent to say  that $v=0$ and
$\epsilon=w$, or that $v\geq 1$ and $\epsilon=w(r-1-v)+1$. In
these cases in $(4)$ one has $R=1$, and therefore from $(1)$ it
follows that $p_a(C)\geq G(r;d,s)$, hence $p_a(C)= G(r;d,s)$. This
proves that $(i)$ implies $(ii)$. Conversely, if property $(ii)$
holds true, then $C$ is a curve with maximal genus $G(r;d,s)$.
Then $C$ is arithmetically Cohen-Macaulay and so $e(C)= max\{i:
h(i)<d\}-1$. A computation proves that this number is just
${\frac{d}{s}}+{\frac{2G(r-1;s)-2-s}{s}}$. This proves that $(ii)$
implies $(i)$.

Moreover, if  $(ii)$ holds, then as before one has $R=1$.
Therefore, using also $(i)$, we obtain $e(C)=(2p_a(C)-2)/d$, and
so $C$ is subcanonical. This proves that $(ii)$ implies $(iii)$.
Conversely, if $p_a(C)= G(r;d,s)$ and $C$ is subcanonical, then
$$
e(C)=(2p_a(C)-2)/d={\frac{d}{s}}+{\frac{2G(r-1;s)-2-s}{s}}+{\frac{2(R-1)}{d}}.
\tag 14
$$
Since $|R|\leq s^3/(r-2)$ and $d>>s$, then previous equality
implies that $R=1$, and so the equality holds in $(5)$. This
proves that $(iii)$ implies $(i)$, and concludes the proof of
Proposition $2$.
\bigskip

\noindent {\bf Remark.} $(i)$ In proving Theorem B, we need the
numerical assumption $d>>s>>t>>u$ only  to use the quoted results
in [CCD], [CCD2] and [D], and to prove $(8)$, $(9)$ and $(13)$. To
this purpose, an elementary computation proves that it suffices
assume that:
$$
d>{\frac{2}{3}}s^4, \quad s>{\frac{2}{3}}t^4, \quad t>408(u+1)^3,
\quad u\geq 2,
$$
or
$$
d>max\left\{{\frac{2}{3}}s^4,\, 12(s+1)^2\right\}, \quad s>t^2-t,
\quad t\geq 2, \quad u=1.
$$

As for the proof of Proposition $1$,  in order to deduce $(i)$ and
$(ii)$ from $(1)$ and [CCD2], it suffices to add to the same
numerical explicit assumption appearing in [CCD2], Theorem $2.2$,
the request that
$$
s_1 > 2(r-l)(l^2+2l+9){\frac{s_2^3}{s_3^2}}s_2\cdot\dots\cdot
s_{l}
$$
(when $l=2$, put $s_{3}=r-3$). Similarly, to prove $(iii)$ we
need the same numerical explicit assumption appearing in [CCD3],
Main Theorem, with the further request that
$$
s_1>{\frac{2s_2^4}{r-2}}.
$$

Also to prove Proposition $2$ we  need the same numerical explicit
assumptions appearing in [CCD3], Main Theorem, with the further
request that
$$
d>{\frac{2s^4}{r-2}},
$$
which enables us to deduce that $R=1$ from $(14)$.

All previous numerical hypotheses are certainly not the best
possible. They only are of the simplest form we were able to
conceive.

\smallskip
$(ii)$ We think it would be interesting to understand  extremal
curves with respect inequality $(9)$ in Lemma $1$. For instance,
when $S$ is a smooth and subcanonical surface, one may prove that
its curves $C$ with extremal speciality are necessarily
subcanonical. In fact if $C$  has speciality $e(C)=
{\frac{d}{s}}+{\frac{2\pi-2-s}{s}}$ then $ p_a(C)\geq
{\frac{d^2}{2s}}+{\frac{d}{2s}}(2\pi-2-s)+1$. On the other hand,
Hodge index Theorem and the assumption on $S$ imply that $
p_a(C)\leq {\frac{d^2}{2s}}+{\frac{d}{2s}}(2\pi-2-s)+1$, and
therefore $ p_a(C)= {\frac{d^2}{2s}}+{\frac{d}{2s}}(2\pi-2-s)+1$.
It follows that $e(C)= {\frac{2p_a(C)-2}{d}}$ and so $C$ is
subcanonical.

Is this conclusion true without assuming $S$ is smooth and
subcanonical ?

\smallskip
$(iii)$ From Proposition $2$ it follows that a curve $C$ with
maximal speciality has arithmetic genus given by the formula $(4)$
with $R=1$. We notice that this condition is not enough to ensure
that the speciality is maximal. In fact, as soon as $\epsilon
=s-1$, in the formula $(4)$ one has $R=1$. Notice also that, in
the particular case $s=r-1$,  Proposition $2$ simply says that for
an integral nondegenerate curve $C\subset \bold P^{r}$ of degree
$d$ one has $e(C)\leq {\frac{d}{r-1}}-{\frac{r+1}{r-1}}$, and the
equality holds if and only if $C$ is a Castelnuovo curve with
$d\equiv 2$ $mod(r-1)$. This is a (we think well known)
consequence of Castelnuovo bound, and it holds for any $d\geq r$.

\smallskip
$(iv)$ In the case there exists some  subcanonical curve verifying
the flag condition $(r-1;s_2,\dots,s_{l})$ and with maximal genus
$G(r-1;s_2,\dots,s_l)$, then one may reformulate property $(iii)$
of Proposition $1$ as follows:
$$
e(r;s_1,\dots,s_l)={\frac{s_1}{s_2}}+ e(r-1;s_2,\dots,s_l)-1,
$$
where $e(r;s_1,\dots,s_l)$ denotes the maximal speciality index
among the curves verifying the flag condition $(r;s_1,\dots,s_l)$.

{

\end
\Refs\widestnumber\key{CCD2}


\ref\key CCD \by L. Chiantini, C. Ciliberto and V. Di Gennaro
\paper The genus of curves in ${\bold P^4}$ verifying certain flag
conditions \jour Manuscripta Math. \vol 88 \yr 1995 \pages 119-134
\endref




\ref\key CCD2 \by L. Chiantini, C. Ciliberto and V. Di Gennaro
\paper On the genus of projective curves verifying certain flag
conditions \jour Boll. U.M.I. \vol (7) 10-B \yr 1996 \pages
701-732
\endref



\ref\key CCD3 \by L.Chiantini, C.Ciliberto and V.Di Gennaro \paper
The genus of projective curves \jour Duke Math. J. \vol 70/2 \yr
1993 \pages 229-245
\endref



\ref\key CD \by C. Ciliberto and V. Di Gennaro \paper Factoriality
of certain threefolds complete intersection in ${\bold P^5}$ with
ordinary double points \jour Commun. Algebra \vol 32 (7) \yr 2004
\pages 2705-2710
\endref



\ref\key D \by  V. Di Gennaro \paper Hierarchical structure of the
family of curves with maximal genus verifying  flag conditions
\jour preprint, math.AG/0504576 \yr 2005
\endref


\ref\key EH \by D. Eisenbud and J. Harris \paper Curves in
Projective Space\jour S\'em. Math. Sup. {\bf 85}, Les Presses du
l'Universit\'e de Montr\'eal, Montr\'eal \yr 1982
\endref




\ref\key GP \by L. Gruson and C. Peskine \paper Genre des courbes
dans l'espace projectif \jour Algebraic Geometry: Proceedings,
Norway, 1977, Lecture Notes in Math., Springer-Verlag, New York
\vol 687 \yr 1978 \pages 31-59
\endref





\endRefs
